\documentclass[12pt,reqno]{amsart}
\usepackage{latexsym,amsmath,amsfonts,amssymb,amsthm}
\theoremstyle{plain}
\newtheorem{Thm}{Theorem}
\newtheorem{Lem}{Lemma}

\theoremstyle{definition}

\theoremstyle{remark}

\def\Z{\mathbb Z}

\def\1{{\bf 1}}
\def\pmod #1{\ ({\rm mod}\ #1)}

\begin{document}
\title{$q$-Analogue of Gauss' Divisibility Theorem}
\author{Hao Pan}
\address{Department of Mathematics, Shanghai Jiaotong University, Shanghai
200240, People's Republic of China} \email{haopan79@yahoo.com.cn}
\begin{abstract} We prove that
$$
\sum_{d\mid n}\mu(d)\frac{(q^{ad};q^{ad})_{n/d}}{(q^d;q^d)_{n/d}}\equiv 0\pmod{[n]_{q^b}},
$$
for any positive integers $n$ and $a$, where $b=(n,a)$.
\end{abstract} \subjclass[2000]{Primary 11A07; Secondary 05A30}
\maketitle

Let $p$ be a prime and $a$ be an integer with $p\nmid a$. The Fermat little theorem asserts that
$$
a^{p-1}\equiv 1\pmod{p}.
$$
For an non-negative integer $n$, define $[n]_q$ by
$$
[n]_q:=\frac{1-q^n}{1-q}=1+q+\cdots+q^{n-1}.
$$
We say $[n]_q$ is the $q$-analogue of the integer $n$ since $\lim_{q\to 1}[n]_q=n$.
If $a$ and $b$ are two positive integers with $a\equiv b\pmod{n}$, then we have
$$
[a]_q=\frac{1-q^a}{1-q}=\frac{1-q^b+q^{b}(1-q^{a-b})}{1-q}\equiv[b]_q\pmod{[n]_q},
$$
where the above congruence is considered in the polynomial ring $\Z[q]$.
Define
$$
(x;q)_n=\begin{cases}
(1-x)(1-xq)\cdots(1-xq^{n-1})&\quad\text{if }n\geq 1,\\
1&\quad\text{if }n=0.
\end{cases}
$$
Then we have a $q$-analogue of Fermat's little theorem:
$$
\frac{(q^a;q^a)_{p-1}}{(q;q)_{p-1}}=\prod_{j=1}^{p-1}[a]_{q^j}=\prod_{j=1}^{p-1}\frac{[aj]_{q}}{[j]_q}\equiv 1\pmod{[p]_q}.
$$
The readers may refer to \cite{Andrews99,Sagan92,PanSun06,Pan07} for more $q$-congruences.

A well-known extension of Fermat's little theorem is Euler's totient
theorem:
$$
a^{\phi(n)}\equiv 1\pmod{n}
$$
for any positive integers $n$ and $a$ with $(a,n)=1$, where $\phi$ is the Euler totient function.
However, another generalization of Fermat's little theorem was found by Gauss in 1863 (cf. \cite[p. 191-193]{SandorCrstici04}):
\begin{equation}
\label{gdt}
\sum_{d\mid n}\mu(d)a^{n/d}\equiv 0\pmod{n}
\end{equation}
for any positive integers $n$ and $a$ (not necessarily co-prime), where $\mu$ is the M\"obius function.
In this short note, we shall give a $q$-analogue of Gauss' divisibility theorem.
\begin{Thm} For any positive integers $n$ and $a$,
$$
\sum_{d\mid n}\mu(d)\frac{(q^{ad};q^{ad})_{n/d}}{(q^d;q^d)_{n/d}}\equiv 0\pmod{[n]_{q^b}},
$$
where $b=(n,a)$ is the greatest common divisor of $n$ and $a$.
\end{Thm}
Let
$$
F(q)=\sum_{d\mid n}\mu(d)\frac{(q^{ad};q^{ad})_{n/d}}{(q^d;q^d)_{n/d}}.
$$
Let $\zeta_m=e^{2\pi\sqrt{-1}/m}$ for each $m\geq 1$. Clearly
$$
[n]_{q^b}=\frac{1-q^{nb}}{1-q^b}=\prod_{\substack{1\leq s\leq nb\\ n\nmid s}}(q-\zeta_{nb}^s).
$$
So it suffices to show that $F(\zeta_{nb}^s)=0$ for each $1\leq s\leq nb$ with $n\nmid s$.
Notice that
$$
\frac{(q^{ad};q^{ad})_{n/d}}{(q^{d};q^{d})_{n/d}}=\prod_{j=1}^{n/d}\frac{1-q^{jad}}{1-q^{jd}}.
$$
Since $\zeta_{nb/d}^{js}=1$ implies $\zeta_{nb/d}^{jsa}=1$, we have
$$
\frac{(q^{ad};q^{ad})_{n/d}}{(q^{d};q^{d})_{n/d}}\bigg|_{q=\zeta_{nb}^s}=0
$$
if there exists $1\leq j\leq n/d$ such that $\zeta_{nb/d}^{jsa}=1$ but $\zeta_{nb/d}^{js}\not=1$.
It follows that $F(\zeta_{nb}^s)=0$ for each $s$ with $b\nmid s$, by noting that
$\zeta_{nb/d}^{(n/d)sa}=\zeta_{b}^{sa}=1$ and
$\zeta_{nb/d}^{(n/d)s}=\zeta_{b}^{s}\not=1$. Now suppose that $b\mid s$ and $t=s/b$.
The following lemma is an easy exercise in elementary number theory.
\begin{Lem}
\label{l1}
$$
|\{1\leq j\leq m:\, jt\equiv 0\pmod{m}\}|=(m,t).
$$
\end{Lem}
By Lemma \ref{l1}
$$
F(\zeta_{n}^{t})=\sum_{\substack{d\mid n\\ (ta,n/d)=(t,n/d)}}\mu(d)\frac{(q^{ad};q^{ad})_{n/d}}{(q^{d};q^{d})_{n/d}}\bigg|_{q=\zeta_{n}^{t}}.
$$
\begin{Lem}
\label{l2}
$$
\prod_{\substack{1\leq j\leq m\\ m\nmid jt}}(1-\zeta_m^{jt})=(m/(m,t))^{(m,t)}.
$$
\end{Lem}
\begin{proof}
\begin{align*}
\prod_{\substack{1\leq j\leq m\\ m\nmid jt}}(1-\zeta_m^{jt})=&\prod_{\substack{1\leq j\leq m\\ j\not\equiv 0\pmod{m/(m,t)}}}(1-\zeta_{m/(m,t)}^{jt/(m,t)})\\
=&\prod_{1\leq j<m/(m,t)}(1-\zeta_{m/(m,t)}^{jt/(m,t)})^{(m,t)}\\
=&(m/(m,t))^{(m,t)}.
\end{align*}
\end{proof}
Now by Lemmas \ref{l1} and \ref{l2},
\begin{align*}
&\sum_{\substack{d\mid n\\ (ta,n/d)=(t,n/d)}}\mu(d)\frac{(q^{ad};q^{ad})_{n/d}}{(q^{d};q^{d})_{n/d}}\bigg|_{q=\zeta_{n}^{t}}\\
=&\sum_{\substack{d\mid n\\ (ta,n/d)=(t,n/d)}}\mu(d)\prod_{\substack{1\leq j\leq n/d\\ (n/d)\nmid jt}}\frac{1-q^{jtad}}{1-q^{jtd}}\bigg|_{q=\zeta_{n}^{t}}
\cdot\prod_{\substack{1\leq j\leq n/d\\ (n/d)\mid jt}}\frac{1-q^{jtad}}{1-q^{jtd}}\bigg|_{q=\zeta_{n}^{t}}\\
=&\sum_{\substack{d\mid n\\ (ta,n/d)=(t,n/d)}}\mu(d)\frac{\prod_{\substack{1\leq j\leq n/d,(n/d)\nmid jta}}(1-\zeta_{n/d}^{jta})}{\prod_{\substack{1\leq j\leq n/d,(n/d)\nmid jt}}(1-\zeta_{n/d}^{jt})}\cdot
\prod_{\substack{1\leq j\leq n/d\\ (n/d)\mid jt}}\bigg(\sum_{k=0}^{a-1}\zeta_{n/d}^{kjt}\bigg)\\
=&\sum_{\substack{d\mid n\\ (ta,n/d)=(t,n/d)}}\mu(d)a^{(n/d,t)}.
\end{align*}
Thus our desired result immediately follows from the next lemma.
\begin{Lem}
\label{l3}
$$
\sum_{\substack{d\mid n\\ (d,ta)=(d,t)}}\mu(n/d)a^{(d,t)}=0
$$
provided that $n\nmid tb$.
\end{Lem}
\begin{proof}
Suppose that $d\mid n$. Clearly $(d,tb)\mid (d,ta)$.
On the other hand,
$$
(d,ta)=((d,n),ta)=(d,(ta,n))\mid (d,t(a,n))=(d,tb).
$$
So we have $(d,ta)=(d,tb)$.
Then
\begin{align*}
\sum_{d\mid n}\mu(n/d)a^{(d,t)}=&
\sum_{u\mid (t,n)}a^u\sum_{\substack{d\mid n, u\mid d\\ (tb/u,d/u)=1}}\mu(n/d)\\
=&\sum_{u\mid (t,n)}a^u\sum_{\substack{d\mid n, u\mid d}}\mu(n/d)\sum_{\substack{v\mid(tb/u,d/u)}}\mu(v)\\
=&\sum_{\substack{u\mid (t,n),\ v\mid (tb/u,n/u)}}a^u\mu(v)\sum_{\substack{d\mid n,\ uv\mid d}}\mu(n/d)\\
=&\sum_{\substack{u\mid (t,n),\ v\mid (tb/u,n/u)}}a^u\mu(v)\sum_{\substack{d\mid (n/uv)}}\mu((n/uv)/d)\\
=&0,
\end{align*}
by noting that $n\not=uv$ since $uv\mid tb$ and $n\nmid tb$.
\end{proof}

\end{document}